\documentclass[a4paper]{article}
\usepackage{hyperref}
\usepackage{cleveref}
\hypersetup{hypertexnames = false, bookmarksdepth = 2, bookmarksopen = true, colorlinks, linkcolor = black, citecolor = black, urlcolor = black, pdfstartview={XYZ null null 1}}

\usepackage{amsfonts}
\usepackage[fleqn,leqno]{amsmath}
\usepackage{amsthm}
\usepackage[maxbibnames=99]{biblatex}
\usepackage{booktabs}
\usepackage{diagbox}
\usepackage{enumitem}
\usepackage{mathtools}
\usepackage{parskip}
\usepackage{pgfplots}
\usepackage{rotating}
\usepackage{siunitx}
\usepackage{thmtools}
\usepackage{tikz-cd}
\usepackage[colorinlistoftodos, textsize = footnotesize]{todonotes}
\usepackage{xparse}
\usepackage{xspace}

\usepackage[T1]{fontenc}
\usepackage{libertine}
\usepackage[libertine]{newtxmath}
\usepackage[scaled=0.83]{beramono}
\usepackage{eucal}
\usepackage{microtype}
\usepackage{gitinfo2}

\newcounter{todocounter}
\DeclareDocumentCommand\addreference{g}{\stepcounter{todocounter}\todo[color = blue!30]{\thetodocounter. Add reference\IfNoValueF{#1}{: #1}}\xspace}
\DeclareDocumentCommand\checkthis{g}{\stepcounter{todocounter}\todo[color = red!50]{\thetodocounter. Check this\IfNoValueF{#1}{: #1}}\xspace}
\DeclareDocumentCommand\fixthis{g}{\stepcounter{todocounter}\todo[color = orange!50]{\thetodocounter. Fix this\IfNoValueF{#1}{: #1}}\xspace}
\DeclareDocumentCommand\expand{g}{\stepcounter{todocounter}\todo[color = green!50]{\thetodocounter. Expand\IfNoValueF{#1}{: #1}}\xspace}

\declaretheoremstyle[
  spaceabove = 3pt,
  spacebelow = 3pt,
]{lecture}
\theoremstyle{lecture}
\newtheorem{theorem}{Theorem}

\newtheorem{definition}[theorem]{Definition}

\newtheorem{proposition}[theorem]{Proposition}
\newtheorem{remark}[theorem]{Remark}

\makeatletter
\def\gitfootnote{\gdef\@thefnmark{}\@footnotetext}
\makeatother

\hypersetup{pagebackref}

\mathchardef\mhyphen="2D
\newcommand\dash{\nobreakdash-\hspace{0pt}}

\DeclareMathOperator\HH{H}

\DeclareMathOperator\moduli{M}
\DeclareMathOperator\Pic{Pic}

\title{Examples violating Golyshev's canonical strip hypotheses}
\author{Pieter Belmans \\ Sergey Galkin \\ Swarnava Mukhopadhyay}
\date{}

\begin{document}
\maketitle

\begin{abstract}
  We give the first examples of smooth Fano and Calabi--Yau varieties violating the (narrow) canonical strip hypothesis, which concerns the location of the roots of Hilbert polynomials of polarised varieties. They are given by moduli spaces of rank~2 bundles with fixed odd-degree determinant on curves of sufficiently high genus, hence our Fano examples have Picard rank~1, index~2, are rational, and have moduli. The hypotheses also fail for several other closely related varieties.
\end{abstract}


\section{The canonical strip hypotheses}
\label{section:hypotheses}
Associated to a polarisation of a smooth projective variety~$X$ we can consider its Hilbert polynomial. The complex roots of this polynomial satisfy a symmetry property induced by Serre duality. In \cite{MR2503098} Golyshev introduced further constraints on these roots: the (narrow) canonical strip hypothesis. The motivation for these restrictions comes from Yau's inequalities on characteristic numbers. At the end of this introduction we give a quick summary of the positive results regarding these hypotheses.

To state (and generalise) the canonical strip hypothesis we will use the following definition.
\begin{definition}
A pair $(X,H)$ of a normal variety and an ample line bundle is said to be \emph{monotone of index~$r$} if
  \begin{equation}
    \mathrm{c}_1(X)=-\mathrm{K}_X\equiv r H,
  \end{equation}
where the symbol $\equiv$ denotes numerical equivalence of divisors.
\end{definition}
The case of~$H=-\mathrm{K}_X$ (resp.~$H=\mathrm{K}_X$) as considered in \cite{MR2503098} for a Fano variety (resp.~variety with~$\mathrm{K}_X$ ample) has index~$1$ (resp.~$-1$). We will also consider polarised Calabi--Yau varieties, for which~$r=0$.

By Serre duality we have that
\begin{equation}
  \chi(nH)=(-1)^{\dim X}\chi(-(r+n)H).
\end{equation}
Hence the roots of the Hilbert polynomial are symmetric around the line~$-r/2$. Golyshev introduced the following further constraints on the real parts of the roots of the Hilbert polynomial.

\begin{definition}
  \label{definition:strip}
  Let~$X$ be a smooth projective variety, and~$H$ an ample line bundle, such that~$(X,H)$ is monotone polarised of index~$r$. Let~$\alpha_1,\ldots,\alpha_{\dim X}$ be the real parts of the roots of the Hilbert polynomial associated to~$H$. Then we say that~$X$ satisfies
  \begin{description}
    \item[(CL)] the \emph{canonical line hypothesis} if
      \begin{equation}
        \alpha_i=r/2,
      \end{equation}
    \item[(NCS)] the \emph{narrow canonical strip hypothesis} if~$r\leq 0$ and
      \begin{equation}
        \alpha_i\in\left[ -r+\frac{r}{\dim X + 1},-\frac{r}{\dim X + 1} \right]
      \end{equation}
      if~$r\geq0$, and
      \begin{equation}
        \alpha_i\in\left[ \frac{-r}{\dim X + 1},-r-\frac{r}{\dim X + 1} \right]
      \end{equation}
      otherwise,
    \item[(CS)] the \emph{canonical strip hypothesis} if~$r\leq 0$ and
      \begin{equation}
        \alpha_i\in[-r,0]
      \end{equation}
      if~$r\geq 0$ and
      \begin{equation}
        \alpha_i\in[0,-r]
      \end{equation}
      otherwise,
  \end{description}
  for all~$i=1,\ldots,\dim X$.
\end{definition}

It is clear that
\begin{equation}
  \mathrm{(CL)}\Rightarrow\mathrm{(NCS)}\Rightarrow\mathrm{(CS)}.
\end{equation}

If~$X$ is a Fano variety, and~$Y\hookrightarrow X$ is a (normal) anticanonical divisor, we can consider the monotone polarised variety~$(Y,-\mathrm{K}_X|_Y)$. By \cite[theorem~4]{MR2503098} we know that if~$(X,-\mathrm{K}_X)$ satisfies~(CS) then~$(Y,-\mathrm{K}_X|_Y)$ satisfies~(CL).

The goal of this paper is to give the first examples of
\begin{enumerate}
  \item Fano varieties which violate the (narrow) canonical strip hypothesis;
  \item embedded Calabi--Yau varieties which violate the canonical line hypothesis.
\end{enumerate}
The question whether such varieties exist was raised by Golyshev in \cite[\S5.A]{MR2503098}. The examples we give are moduli spaces~$\moduli_C(2,\mathcal{L})$ of vector bundles of rank~2 with fixed determinant~$\mathcal{L}$ of odd degree on a curve~$C$ of genus~$g\gg 2$.

\begin{theorem}
  \label{theorem:main-theorem}
  We have the following examples violating the (narrow) canonical strip hypothesis.

  \begin{itemize}
    \item Let~$g\geq 8$, then~$\moduli_C(2,\mathcal{L})$ does not satisfy the narrow canonical strip hypothesis.
    \item Let~$g\geq 10$, then~$\moduli_C(2,\mathcal{L})$ does not satisfy the canonical strip hypothesis.
    \item Let~$g\geq 11$ then an anticanonical Calabi--Yau hypersurface inside~$\moduli_C(2,\mathcal{L})$ does not satisfy the canonical line hypothesis\footnote{Hence for~$g=10$ we have that~$\moduli_C(2,\mathcal{L})$ violates the canonical strip hypothesis, yet an anticanonical Calabi--Yau hypersurface still satisfies the embedded canonical line hypothesis. See also \cref{table:values} for more information.}.
  \end{itemize}
\end{theorem}
Observe that there exist smooth anticanonical hypersurfaces, by the very ampleness of~$\Theta$ \cite{MR1739370} and the Bertini theorem.

In \cref{section:examples} we give the proof of this theorem, and discuss related constructions, giving more families of examples vilating Golyshev's hypotheses. Before we do this we give an overview of the positive results in the literature. In \cite{MR2503098} Golyshev explains how
\begin{enumerate}
  \item the canonical line hypothesis holds for smooth projective curves (with the elliptic curve being embedded in~$\mathbb{P}^2$);
  \item the narrow canonical strip hypothesis holds for del Pezzo surfaces and surfaces of general type, and the canonical line hypothesis holds for embedded K3~surfaces;
  \item the narrow canonical strip hypothesis holds for Fano 3-folds and minimal threefolds of general type.
\end{enumerate}
Moreover it is explained how all Grassmannians (not just projective spaces) satisfy the narrow canonical strip hypothesis.

In \cite{0904.2470v1} Manivel shows that for~$G$ a simple affine algebraic group and~$P$ a maximal parabolic subgroup
\begin{enumerate}
  \item $G/P$ satisfies the tight\footnote{A strengthening of the narrow canonical strip hypothesis for Fano varieties involving the index~$\imath_X$ of~$X$, i.e. with the notation of \cref{definition:strip} one asks for~$\alpha_i\in[-1+1/{\imath_X}\leq -1/{\imath_X}$, when~$H=-\mathrm{K}_X$.} strip hypothesis;
  \item Fano complete intersections in~$G/P$ satisfy the tight canonical strip hypothesis;
  \item general type complete intersections in~$G/P$ satisfy the canonical line hypothesis;
  \item Calabi--Yau complete intersections in~$G/P$ satisfy the canonical line hypothesis.
\end{enumerate}

Miyaoka's celebrated pseudo-effectivity theorem \cite{MR946247} implies that the embedded canonical line hypothesis holds for smooth Calabi--Yau threefolds\footnote{As well as for threefolds with numerically trivial canonical bundle, and terminal Gorenstein singularities that admit crepant resolution.}.

Another case that can be checked is that of smooth toric Fano~$n$\dash folds, for~$n=4,\ldots,7$. By \cite[proposition~9.4.3]{MR2810322} we have that the Hilbert polynomial for the anticanonical bundle is the Ehrhart polynomial of the moment polytope. In \cite{github-toric} we have computed these Ehrhart polynomials, based on the classification of the toric Fano polytopes up to dimension~7. It turns out there are no examples violating the canonical strip hypothesis. In other words, we can add the following proposition to the list of positive examples.

\begin{proposition}
  \label{proposition:no-toric-counterexamples}
  Let~$X$ be a smooth toric Fano variety of dimension at most~7. Then~$X$ satisfies the canonical strip hypothesis\footnote{The narrow canonical strip hypothesis is violated starting in dimension~4.}.
\end{proposition}

The maximal value~$m_d$ of the real parts of the roots of the Hilbert polynomials for smooth toric Fano varieties of dimension~$d$ is given as
\begin{equation}
  \begin{aligned}
    m_2&=-0.3333333333\ldots \\
    m_3&=-0.2500000000\ldots \\
    m_4&=-0.1394448724\ldots \\
    m_5&=-0.0868988066\ldots \\
    m_6&=-0.0566708554\ldots \\
    m_7&=-0.0354049073\ldots
  \end{aligned}
\end{equation}

\paragraph{Acknowledgements}
The first and third author were supported by the Max Planck Institute for Mathematics in Bonn. The second author was partially supported by the Hausdorff Center for Mathematics during the trimester program ``Periods in Number Theory, Algebraic Geometry and Physics'' and by the Laboratory of Mirror Symmetry NRU HSE, RF Government grant, ag.~N.~14.641.31.0001.

\section{Examples violating the hypotheses}
\label{section:examples}
An interesting class of Fano varieties is given by moduli spaces of vector bundles on a curve. We will restrict ourselves to the case of rank~2. Let~$\mathcal{L}$ be a line bundle of odd degree on a smooth projective curve~$C$ of genus~$g$. Then the moduli space~$\moduli_C(2,\mathcal{L})$ of rank~2 bundles with determinant~$\mathcal{L}$ is a smooth projective variety of dimension~$3g-3$, of rank~1 and index~2, i.e.~$\Pic\moduli_C(2,\mathcal{L})\cong\mathbb{Z}\Theta$, and~$\omega_{\moduli_C(2,\mathcal{L})}\cong\Theta^{\otimes-2}$, see \cite{MR999313}.

To compute the Hilbert polynomial we can use the celebrated Verlinde formula, which gives an expression for~$\dim_k\HH^0(\moduli_C(2,\mathcal{L}),\Theta^{\otimes k})$, see \cite{MR1397056,MR1360519} for a survey. It reads
\begin{equation}
  \dim_k\HH^0(\moduli_C(2,\mathcal{L}),\Theta^{\otimes k})
  =
  (k+1)^{g-1}\sum_{j=1}^{2k+1}\frac{(-1)^{j-1}}{\sin^{2g-2}\frac{j\pi}{2k+2}}.
\end{equation}
Rather than this version of the Verlinde formula we will use an alternative form, taken from \cite{MR1360519}. Namely item~(x) in theorem~1 of op.~cit.~gives the formula
\begin{equation}
  \dim_k\HH^0(\moduli_C(2,\mathcal{L}),\Theta^{\otimes k})
  =
  \frac{2^g\det M_{r,s}}{\prod_{j=1}^g(2j)!}
\end{equation}
where the matrix~$(M_{r,s})_{r,s=0,\ldots,g-1}$ is given by
\begin{equation}
  M_{r,s}=
  \begin{cases}
    1 & r=0 \\
    (k+1+r)^{2s+2} - (k+1-r)^{2s+2} & r\geq 0
  \end{cases}.
\end{equation}
The benefit of using this expression is that it can be computed in an exact fashion in computer algebra.

Using this formula one computes the first~$3g$ coefficients of the Hilbert series, and from this we can obtain the Hilbert polynomial of~$\moduli_C(2,\mathcal{L})$ with respect to~$\Theta$, i.e.~we consider the monotone polarisation given by~$H=\Theta$ for~$\moduli_C(2,\mathcal{L})$. Two implementations of the computations (one in Pari/GP, another in Sage) can be found at \cite{github-moduli}. The implementation computes the maximum of the real parts of the complex roots of the Hilbert polynomial, so we are interested in knowing when these are negative, but close to~0, or positive. From these computations we get \cref{theorem:main-theorem} as in the introduction.

\begin{sidewaysfigure}
  \input{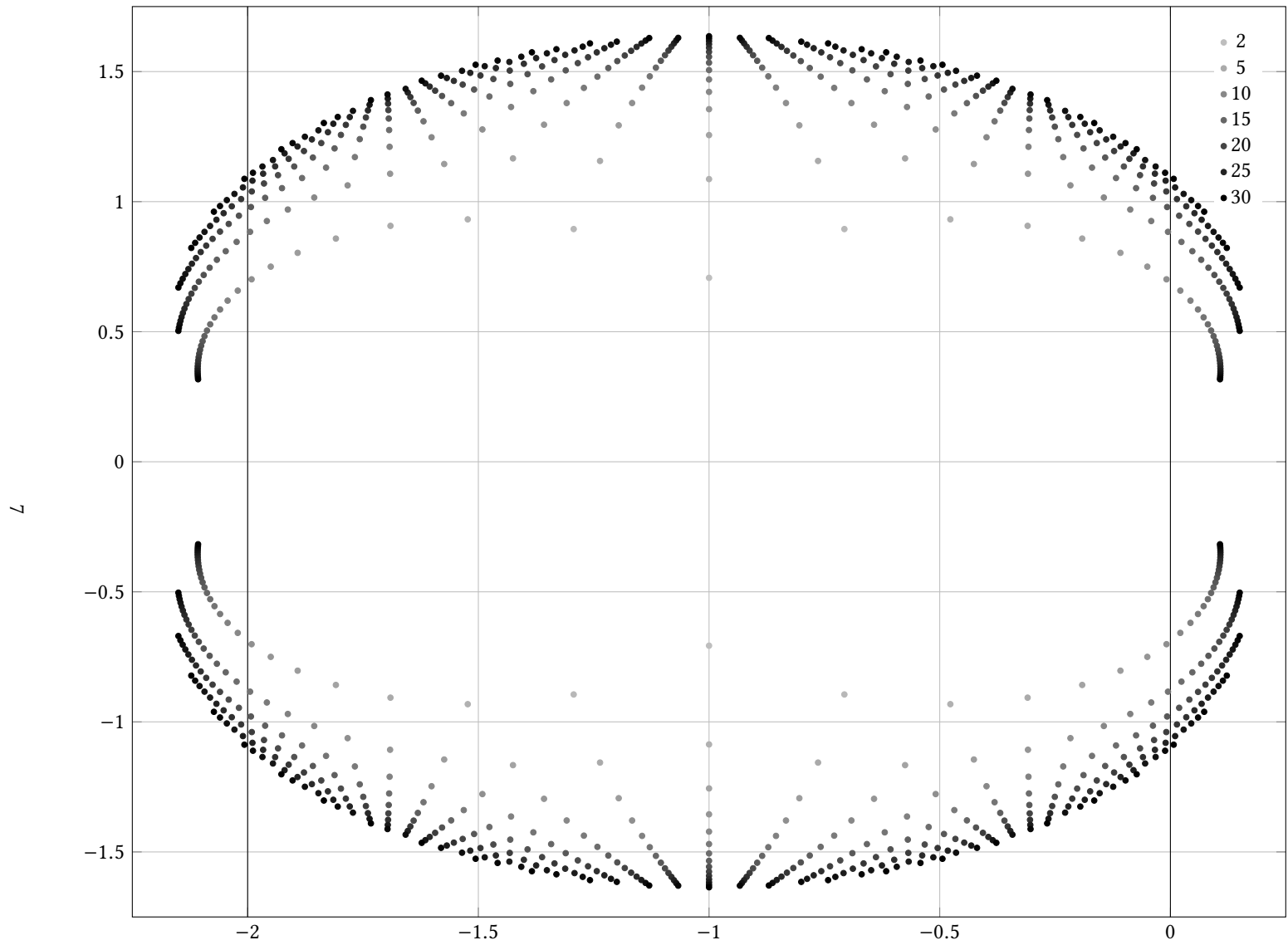}
  \caption{Complex roots of Hilbert polynomials of $\moduli_C(2,\mathcal{L})$, for $g=2,\ldots,30$}
  \label{figure:plot}
\end{sidewaysfigure}

\begin{remark}
  The values in the column labeled~$\moduli_C(2,\mathcal{L})$ in \cref{table:values} suggest an interesting convergence behaviour for the maximum of the real part of the complex roots of the Hilbert polynomial. More generally one can compute that the collection of all roots of the Hilbert polynomial seems to exhibit a pattern where the limiting behaviour involves the complex hull of the roots for increasing genera. A visualisation of this is given in \cref{figure:plot}. In the picture we have omitted the root at~$t=-1$, which in all the examples we computed is of multiplicity~$g-1$, but we have no proof of this. We suggest these questions for future work.
\end{remark}

\paragraph{Related constructions}
Besides an anticanonical Calabi--Yau hypersurface constructed out of~$\moduli_C(2,\mathcal{L})$ there are other Fano and Calabi--Yau varieties we can construct out of it. These are
\begin{description}[leftmargin=1cm,labelindent=.5cm]
  \item[Fano\textsuperscript{1}] the~$3g-4$\dash dimensional Fano variety given by a linear section;
  \item[Fano\textsuperscript{2}] the~$3g-3$\dash dimensional Fano variety given by a double cover branched in~$2\Theta$;
  \item[CY\textsuperscript{2}] the~$3g-5$\dash dimensional Calabi--Yau variety given by a linear section of codimension~2;
  \item[CY\textsuperscript{3}] the~$3g-3$\dash dimensional Calabi--Yau variety given by a double cover branched in~$4\Theta$;
  \item[CY\textsuperscript{4}] the~$3g-3$\dash dimensional Calabi--Yau variety given by the cone over the embedding given by~$\Theta$, intersected with a cubic hypersurface;
  \item[CY\textsuperscript{5}] the~$3g-3$\dash dimensional Calabi--Yau variety given by the join with a line intersected with two quadric hypersurfaces;
  \item[CY\textsuperscript{6}] the~$3g-3$\dash dimensional Calabi--Yau variety given by a smoothing of a linear section of a join with an elliptic curve of degree 1.
\end{description}
For all of these the canonical strip (resp.~line) hypothesis eventually fails, as checked in \cite{github-moduli}. In \cref{table:values} we have collected the maximum over the real parts of the complex roots of the Hilbert polynomial, where the columns are labelled as in this remark. The Calabi--Yau variety denoted~CY\textsuperscript{1} is the anticanonical section of~$\moduli_C(2,\mathcal{L})$ as considered in \cref{theorem:main-theorem}.

We have not found counterexamples with ample canonical bundle: the canonical line hypothesis was satisfied for all constructions we considered.

\begin{sidewaystable}
  \centering
    \small
    \begin{tabular}{cS[table-format=1.10]S[table-format=1.10]|S[table-format=1.10]|S[table-format=1.10]S[table-format=1.10]S[table-format=1.10]S[table-format=1.10]S[table-format=1.10]S[table-format=1.10]}
    & \multicolumn{3}{c|}{\textbf{Fano}} & \multicolumn{6}{c}{\textbf{Calabi--Yau}} \\
    \toprule
    $g$ & {Fano\textsuperscript{1}} & {Fano\textsuperscript{2}} & {$\moduli_C(2,\mathcal{L})$} & CY\textsuperscript{1} & {CY\textsuperscript{2}} & {CY\textsuperscript{3}} & {CY\textsuperscript{4}} & {CY\textsuperscript{5}} & {CY\textsuperscript{6}} \\
    \midrule
2  & -0.5          & -0.5          & -1            & 0            & 0            & 0            & 0            & 0            & 0 \\
3  & -0.5          & -0.5          & -0.7066405395 & 0            & 0            & 0            & 0            & 0            & 0 \\
4  & -0.5          & -0.5          & -0.4770019488 & 0            & 0            & 0            & 0            & 0            & 0 \\
5  & -0.2890507098 & -0.3131727064 & -0.3094989272 & 0            & 0            & 0            & 0            & 0            & 0 \\
6  & -0.1792056326 & -0.2063905610 & -0.1911961780 & 0            & 0            & 0            & 0            & 0            & 0 \\
7  & -0.1047144340 & -0.1119844025 & -0.1083536780 & 0            & 0            & 0            & 0            & 0            & 0 \\
8  & -0.0500408825 & -0.0499879643 & -0.0500409722 & 0            & 0            & 0            & 0            & 0            & 0 \\
9  & -0.0088875090 & -0.0081356074 & -0.0085094225 & 0            & 0            & 0            & 0            & 0            & 0 \\
10 & 0.0213534238  & 0.0216201879  & 0.0214869361  & 0            & 0.0379539521 & 0.0381695630 & 0.0382767453 & 0.0383835172 & 0.0380619666 \\
11 & 0.0434392549  & 0.0434699963  & 0.0434546003  & 0.0614369091 & 0            & 0            & 0            & 0            & 0 \\
12 & 0.0597507064  & 0.0597412231  & 0.0597459652  & 0.0399471632 & 0.0731794361 & 0.0731745236 & 0.0731720669 & 0.0731696099 & 0.0731769800 \\
13 & 0.0719600677  & 0.0719550941  & 0.0719575810  & 0.0801077393 & 0.0675677156 & 0.0675620782 & 0.0675592588 & 0.0675564391 & 0.0675648971 \\
14 & 0.0811899396  & 0.0811890603  & 0.0811894999  & 0.0819305430 & 0.0845735173 & 0.0845730490 & 0.0845728148 & 0.0845725807 & 0.0845732831 \\
15 & 0.0882121052  & 0.0882121423  & 0.0882121238  & 0.0879245273 & 0.0907743344 & 0.0907742826 & 0.0907742567 & 0.0907742308 & 0.0907743085 \\
16 & 0.0935738073  & 0.0935738646  & 0.0935738359  & 0.0965258891 & 0.0911200604 & 0.0911201236 & 0.0911201551 & 0.0911201867 & 0.0911200920 \\
17 & 0.0976711255  & 0.0976711387  & 0.0976711321  & 0.0946222779 & 0.1003675084 & 0.1003675211 & 0.1003675275 & 0.1003675339 & 0.1003675148 \\
18 & 0.1007949361  & 0.1007949368  & 0.1007949365  & 0.1029737199 & 0.0981849016 & 0.0981849019 & 0.0981849020 & 0.0981849022 & 0.0981849018 \\
19 & 0.1058859249  & 0.1058863358  & 0.1058861304  & 0.1051019381 & 0.1070028490 & 0.1070031860 & 0.1070033546 & 0.1070035231 & 0.1070030175 \\
20 & 0.1146393484  & 0.1146393524  & 0.1146393504  & 0.1144075091 & 0.1150500957 & 0.1150501150 & 0.1150501247 & 0.1150501344 & 0.1150501053 \\
21 & 0.1218850498  & 0.1218850164  & 0.1218850331  & 0.1225735595 & 0.1211992074 & 0.1211991712 & 0.1211991532 & 0.1211991351 & 0.1211991893 \\
22 & 0.1278911325  & 0.1278911199  & 0.1278911262  & 0.1272829480 & 0.1284698807 & 0.1284698686 & 0.1284698625 & 0.1284698564 & 0.1284698747 \\
23 & 0.1328722016  & 0.1328721997  & 0.1328722006  & 0.1332346075 & 0.1324975586 & 0.1324975566 & 0.1324975556 & 0.1324975546 & 0.1324975576 \\
24 & 0.1370012165  & 0.1370012167  & 0.1370012166  & 0.1368306124 & 0.1371716714 & 0.1371716716 & 0.1371716717 & 0.1371716719 & 0.1371716715 \\
25 & 0.1404184745  & 0.1404184747  & 0.1404184746  & 0.1404629798 & 0.1403761630 & 0.1403761632 & 0.1403761633 & 0.1403761634 & 0.1403761631 \\
    \midrule
    $\dim$ & {$3g-4$} & {$3g-3$} & {$3g-3$} & {$3g-4$} & {$3g-5$} & {$3g-3$} & {$3g-3$} & {$3g-3$} & {$3g-3$} \\
    \bottomrule
  \end{tabular}
  \caption{Maximum value of real parts of complex roots of Hilbert polynomial}
  \label{table:values}
\end{sidewaystable}

\clearpage

\printbibliography

\textsc{Max-Planck-Institut fur Mathematik, Bonn, Germany} \\
\texttt{pbelmans@mpim-bonn.mpg.de} \\
\texttt{swarnava@mpim-bonn.mpg.de}

\textsc{National Research University Higher School of Economics, Russian Federation}
\texttt{sergey.galkin@phystech.edu}

\end{document}